# Combining heuristics and Exact Algorithms: A Review

Hengameh Fakhravar
Old Dominion University
Hfakh001@odu.edu

*Abstract*

Several different ways exist for approaching hard optimization problems. Mathematical programming techniques, including (integer) linear programming based methods, and metaheuristic approaches are two highly successful streams for combinatorial problems. These two have been established by different communities more or less in isolation from each other. Only over several years ago a larger number of researchers recognized the advantages and huge potentials of building hybrids of mathematical programming methods and metaheuristics. In fact, many problems can be practically solved much better by exploiting synergies between these different approaches than by "pure" traditional algorithms. The crucial issue is how mathematical programming methods and metaheuristics should be combined for achieving those benefits. This paper surveys existing techniques for such combinations and provide some examples of using them for vehicle routing problem.

*Index Terms* matheuristics, optimization-based heuristics, VRP, survey,

## INTRODUCTION

Many problems arising in areas such as scheduling and production planning, location and distribution management, Internet routing or bioinformatics are combinatorial optimization problems (COPs). COPs are intriguing because they are often easy to state but often very difficult to solve, which is captured by the fact that many of them are NP-hard [1]. This difficulty and, at the same time, their enormous practical importance, have led to a large number of solution techniques for them. The available solution techniques can be classified as being either exact or approximate algorithms. Exact algorithms are guaranteed to find an optimal solution and prove its optimality for every finite size instance of a COP within an instance-dependent, finite run-time, or prove that no feasible solution exists. If optimal solutions cannot be computed efficiently in practice, it is usual to trade the guarantee of optimality for efficiency. In other words, the guarantee of finding optimal solutions is sacrificed for the sake of getting very good solutions in reasonably short time by using approximate algorithms.

Two solution method classes that have significant success are integer programming (IP) methods as an exact approach, and stochastic local search (SLS) algorithms [2] as an approximate approach. IP methods rely on the characteristic of the decision variables being integers. Some well-known IP methods are branch-and-bound, branch-and-cut, branch-and-price, and dynamic programming. Important advantages of exact methods for IP are that (i) proven optimal solutions can be obtained if the algorithm succeeds, (ii) valuable information on upper/lower bounds on the optimal solution are obtained even if the algorithm is stopped before completion (IP methods can become approximate if we define a criterion of stopping them before solving the problem), and (iii) IP methods allow to provably prune parts of the search space in which optimal solutions cannot be located. A more practical advantage of IP methods is that research code such as Minto [3] or GLPK [4], or powerful, general-purpose commercial tools such as CPLEX [5] or Xpress-MP [6] are available. However, despite the known successes, exact methods have a number of disadvantages. Firstly, for many problems the size of the instances that are practically solvable is rather limited and, even if an application is feasible, the variance of the computation times is typically very large when applied to different instances of a same size. Secondly, the memory consumption of exact algorithms can be very large and lead to the early abortion of a program. Thirdly, for many COPs the best performing algorithms are problem specific and they require large development times by experts in integer programming. Finally, high performing exact algorithms for one problem are often difficult to extend if some details of the problem formulation change. The state-of-the art for exact algorithms is that for some NP-hard problems very large instances can be solved fast, while for other problems even small size instances are out of reach.

SLS is probably the most successful class of approximate algorithms. When applied to hard COPs, local search yields high-quality solutions by iteratively applying small modifications (local moves) to a solution in the hope of finding a better one. Embedded into higher-level guidance mechanisms, which are called (general-purpose) SLS methods [2] or, more commonly, metaheuristics, this approach has been shown to be very successful in achieving near-optimal (and often optimal) solutions to a number of difficult problems[2] [7]. Examples of well-known general-purpose SLS methods (or metaheuristics) are simulated annealing, tabu search, memetic algorithms, ant colony optimization or iterated local search[8]. Advantages of SLS algorithms are that (i) they are the best performing algorithms available for a variety of problems, (ii) they can examine a huge number of possible solutions in short computation time, (iii) they are often more easily adapted to slight variants of problems and are therefore more flexible, and (iv) they are typically easier to understand and implement by the common user than exact methods. However, local search based algorithms have several disadvantages. Firstly, they cannot prove optimality and typically do not give bounds on the quality of the solutions they return. Secondly, they typically cannot provably reduce the search space. Thirdly, they do not have well defined stopping criteria (this is particularly true for metaheuristics). Finally, local search methods often have problems with highly constrained problems where feasible areas of the solution space are disconnected. Another problem that occurs in practice is that, with very few exceptions [9], there are no efficient general-purpose local search solvers available. Hence, although one can typically develop an SLS algorithms of reasonable performance rather quickly, many applications of SLS algorithms can require considerable development and implementation efforts if very high performance is required.

It should be clear by now that IP and SLS approaches have their particular advantages and disadvantages and can be seen as complementary. Therefore, it appears to be a good idea to try to combine these two distinct techniques into more powerful algorithms.

When considering optimization approaches that combine aspects from metaheuristics with mathematical programming techniques, the resulting hybrid system may either be of exact or heuristic nature. Exact approaches are guaranteed to yield proven optimal solutions when they are given enough computation time. In contrast, heuristics only aim at finding reasonably good approximate solutions usually in a more restricted time; performance guarantees are typically not provided. Most of the existing hybrid approaches are of heuristic nature, and mathematical programming techniques are used to boost the performance of a metaheuristic. Exploiting solutions to exactly solvable relaxations of the original problem, or searching large neighborhoods by means of mathematical programming techniques are examples for such approaches.

On the other hand, there are also several highly successful ways to exploit metaheuristic strategies for enhancing the performance of mathematical programming techniques, and often these methods retain their exactness.

In the first section, we will continue with a structural classifications of strategies for combining metaheuristics and exact optimization techniques. In the second section, we discuss the matheuristic approaches for the routing problems. The last section is devoted for general discussion and conclusion.

### STRUCTURAL MODELS FOR COMBINING METAHEURISTICS WITH EXACT APPROACHES

The available techniques for COPs can roughly be classified into two main categories: exact and heuristic methods. Exact algorithms are guaranteed to find an optimal solution and to prove its optimality for every instance of a COP. The run-time, however, often increases dramatically with the instance size, and often only small or moderately-sized instances can be practically solved to provable optimality. In this case, the only possibility for larger instances is to trade optimality for run-time, yielding heuristic algorithms. In other words, the guarantee of finding optimal solutions is sacrificed for the sake of getting good solutions in a limited time. Two independent heterogeneous streams, coming from very different scientific communities, had significant success in solving COPs:

– Integer Programming (IP) as an exact approach, coming from the operations research community and based on the concepts of linear programming[10]. Among the exact methods are branch-and-bound (B&B), dynamic programming, Lagrangian relaxation based methods, and linear and integer programming-based methods, such as branch-and-cut, branch-and-price, and branch-and-cut and-price [11].

– Local search with various extensions and independently developed variants, in the following called metaheuristics, as a heuristic approach. Metaheuristics include, among others, simulated annealing [12], tabu search [13], iterated local search [14], variable neighborhood search [15], and various population-based models such as evolutionary algorithms [16], scatter search [17], memetic algorithms [18], and various estimation of distribution algorithms [19].

In [20] authors present a more general classification of existing approaches combining exact and metaheuristic algorithms for combinatorial optimization in which the following two main categories are distinguished:

– *Collaborative Combinations*: By collaboration it means that the algorithms exchange information, but are not part of each other. Exact and heuristic algorithms may be executed sequentially, intertwined or in parallel.

– *Integrative Combinations*: By integration it means that one technique is a subordinate embedded component of another technique. Thus, there is a distinguished master algorithm, which can be either an exact or a metaheuristic algorithm, and at least one integrated slave.

[21] present a similar classification of hybrid algorithms, further including constraint programming. The authors discern a decomposition scheme corresponding to the integrative combinations and a multiple search scheme corresponding to collaborative combinations.

In another classification, [22] classifies heuristics approaches to 4 category and then tried to show how we can use mathematical programming in each.

– *Construction heuristics*: start from ''scratch'' and proceed through a set of steps, each of which adds a component to the solution until a complete (feasible) solution is generated. We also label such methods decomposition approaches since they effectively decompose a larger problem into a series of sequentially executed sub problems.

– *Improvement heuristics*: start with a feasible solution and iteratively execute solution improving steps until some termination condition is met.

– *Relaxation-based heuristics*: It is very often the case that while a problem may be very difficult, a certain relaxation to that problem may be efficiently solvable. The solution to a relaxation generates a bound on the value of a problem's optimal solution. As such relaxations are often employed in exact mathematical programming approaches. Additionally, they can often serve as a basis for effective heuristics. Two general approaches are used. In one, the solution to a relaxation is modified to generate a feasible solution to the problem of interest. Probably the prototypical approach of this type involves rounding of the solution to a linear programming relaxation of an integer program. The second class of relaxation based approaches makes use of the dual information provided by the solution to the relaxation in a subsequently executing heuristic.

– *Using mathematical programming algorithms to generate approximate solutions*: An exact optimization algorithm terminates with an optimal solution and a proof of optimality. In many cases, a significant portion of the total solution time is spent proving that a solution found (quickly) is optimal. Another common scenario is that a large amount of computing time is spent going from a ''near optimal'' solution to an optimal one. With this motivation, in many practical settings, exact mathematical programming algorithms are modified to generate a very good, but not necessarily optimal, solutions. Approaches of this class are based on the idea of solving the mathematical programming formulation in a 'relaxed' way, i.e., by relaxing some attributes of the exact solution approach that increase solution time significantly. Methodologies that fall in this class are, for example, the premature stopping of a branch-and-bound algorithm, heuristic variable fixing, and rounding of the relaxed solution. Also, the branch-and-price/column generation-based approaches belong to this class.

The other survey on matheuristics is the one done by [23]. The proposed classification is different from the one adopted in [22]. The following classes are introduced:

1. *set-covering/partitioning-based approaches*, corresponding to the class of branch and- price/column generation-based approaches;
2. *Local branching approaches*, based on the local branching scheme proposed in [24].
3. *Decomposition approaches*, which coincides with the first class defined in [22].

## MATHEURISTICS FOR VEHICLE ROUTING PROBLEM: A REVIEW

[25] classify matheuristics for vehicle routing problems into three classes, which we state verbatim:
1. *Decomposition approaches*. In general, in a decomposition approach the problem is divided into smaller and simpler sub-problems and a specific solution method is applied to each sub-problem. In matheuristics, some or all these sub-problems are solved through mathematical programming models to optimality or sub-optimality.
2. *Improvement heuristics*. Matheuristics belonging to this class use mathematical programming models to improve a solution found by a different heuristic approach. They are very common as they can be applied whatever heuristic is used to obtain a solution that the mathematical programming model aims at improving.
3. *Branch-and-price/column generation-based approaches*. Branch-and-price algorithms have been widely and successfully used for the solution of routing problems. Such algorithms make use of a set partitioning formulation, where a binary or integer variable is associated with each possible route (column). Due to the exponential number of variables, the solution of the linear relaxation of the formulation is performed through column generation. In the branch-and-price/column generation-based matheuristics the exact method is modified to speed up the convergence, thus losing the guarantee of optimality. For example, the column generation phase is stopped prematurely.

The following three sections are devoted to the description of these three classes.

### I. Decomposition Approaches

Traditionally, heuristic methods, and metaheuristics in particular, have been primal-only methods. They are usually quite effective in solving the given problem instances, and they terminate providing the best feasible solution found during the allotted computation time. However, disregarding dual information implies some obvious drawbacks, first of all not knowing the quality of the proposed solution, but also having possibly found an optimal solution at the beginning of the search and having wasted CPU time ever since, having searched a big search space that could have been much reduced, or having disregarded important information that could have been very effective for constructing good solutions. Dual information is also tightly connected with the possibility of obtaining good lower bounds (making reference, here and forward, to minimization problems), another element which is not a structural part of current metaheuristics. On the contrary, most mathematical programming literature dedicated to exact methods is strongly based on these elements for achieving the obtained results. There is nothing, though, that limits the effectiveness of dual/bounding procedures to exact methods. There are in fact wide research possibilities both in determining how to convert originally exact methods into efficient heuristics and in designing new, intrinsically heuristic techniques, which include dual information.

There are many ways in which bounds can be derived, one of the most effective of these is the use of decomposition techniques [25], [26]. These are techniques primarily meant to exploit the possibility of identifying a sub-problem in the problem to solve and to decompose the whole problem in a master problem and a sub-problem, which communicate via dual or dual-related information. The sub-problems are handled and solved independently. Finally, a feasible solution for the original problem is obtained from the solutions to the sub-problems. In matheuristics, one or all the sub-problems are solved through the exact solution of a mathematical programming formulation. There are three basic decomposition techniques: *Lagrangean relaxation*, *Dantzig- Wolfe decomposition*, and *Benders decomposition.* The popularity of these techniques derives both from their effectiveness in providing efficient bounds and from the observation that many real-world problems lead themselves to a decomposition.

Unfortunately, despite their prolonged presence in the optimization literature, there is as yet no clear-cut recipe for determining which problems should be solved with decompositions and which are better solved by other means. Clearly, decomposition techniques are foremost candidates for problems which are inherently structured as a master and different sub-problems, but it is at times possible to effectively decompose the formulation of a problem which does not show such structure and enjoy advantages. Examples from the literature of effective usage of decomposition techniques (mainly Lagrangean) on single-structure problems include, e.g., set covering [27], [28],[29] set partitioning [30], [31], [32],[33] and crew scheduling [34], [35], [36], [37][38]. This is also the case for the vehicle routing problems (VRPs), the inventory routing problems (IRPs), the production routing problems (PRPs) and the location routing problems (LRPs). Different matheuristics have been proposed for the solution of these problems belonging to the class of decomposition approaches. Routing problems typically involve the following two basic decisions (in addition to further decisions related to the particular application): the clustering of customers which are assigned to each vehicle and the sequencing of customers in vehicle routes. This feature makes it natural the use of a decomposition approach of the kind cluster first-route second, i.e., an approach where first the assignment of customers to vehicles is made and then the decision on how to route the customers assigned to each vehicle is taken.

One of the most used approaches for routing problem is *cluster first-route second* approach. [39],[40]

- *Cluster first-route second approaches*

The basic idea of the cluster first-route second approaches is to divide the two main decisions that characterize routing problems, i.e., the assignment of customers to vehicles and the sequencing of the customers visited by each route. Cluster first-route second is one of the first heuristic approaches proposed for the solution of the classical VRP. In the VRP, we are given a set of customers with an associated demand and a fleet of capacitated vehicles. The problem is to find a set of vehicle routes serving these customers such that the demands are satisfied, each customer is visited only once and the vehicle capacity is never exceeded.

A matheuristic based on a cluster first-route second approach to solve the VRP is motivated by the fact that the clustering of customers can be handled through the solution of a MILP. The routing of customers inside each route can instead be managed by adopting any heuristic available for the solution of the Traveling Salesman Problem (TSP).

The first authors who proposed a cluster first-route second matheuristic for a routing problem, and specifically for the VRP, are [36]. In the first phase of the algorithm so-called seed customers are chosen heuristically, and an assignment problem is solved to optimality to assign the remaining customers to the seed customers. Each seed customer identifies a cluster of customers. Then, routes are built by solving a TSP on each cluster. This approach can be applied to a wide variety of routing problems. The scheme has later been extended to solve the VRPTW in [37]. [38]propose a decomposition approach for the VRP which is similar to the one proposed in [36]. The algorithm is based on the formulation of the routing problem as a capacitated concentrator location problem (CCLP). The idea is to identify seed points, to estimate the cost of assigning each customer to each seed point and then solve a CCLP to determine the clustering of customers. Once the clusters are obtained, a TSP is solved on each cluster. The authors apply the algorithm to the VRP showing that the heuristic performs well on both problems and often outperforms previous heuristics proposed in the literature. A similar approach is applied by the same authors to the VRPTW in [41],[42]

## II. Improvement Heuristics

Matheuristics belonging to the class of improvement heuristics combine a heuristic with the exact solution of a MILP model that aims at improving the solution obtained by applying the heuristic. Different ways to combine the heuristic procedure and the solution of a MILP model have been developed. This combination can go two-ways, either using MILP to improve or design metaheuristics or using metaheuristics for improving known MILP techniques, even though the first of these two directions is by far more studied.

When using MILP embedded into metaheuristics, the main possibility appears to be improving local search [44],[45]. A seminal work in this direction is local branching [46],[47], where MILP is used to define a suitable neighborhood to be explored exactly by a MILP solver. Essentially, only a number of decision variables is left free and the neighborhood is composed by all possible value combination of these free variables.

The idea of an exact exploration of a possibly exponential size neighborhood is at the heart of several other approaches. One of the best known is possibly *Very Large Neighborhood Search* (VLNS) [48],[49]. This method can be applied when it is possible to define the neighborhood exploration as a combinatorial optimization problem itself. In this case it could be possible to solve it efficiently, and it becomes possible the full exploration of exponential neighborhoods. Complementary to this last is the *corridor* method [50],[51], [52] where a would-be large exponential neighborhood is kept of manageable size by adding exogenous constraint to the problem formulation, so that the feasible region is reduced to a "corridor" around the current solution.

Several other methods build around the idea of solving via MILP the neighborhood exploration problem, they differ in the way the neighborhood is defined. For example, an unconventional way of defining it is proposed in the '*dynasearch*' method [53],[54], where the neighborhood is defined by the series of moves which can be performed at each iteration, and dynamic programming is used to find the best sequence of simple moves to use at each iteration.

However, MILP contributed to metaheuristics also along two other opposite lines: improving the effectiveness of well-established metaheuristics and providing the structural basis for designing new metaheuristics. As for the first line, MILP hybrids are reported for most known metaheuristics: tabu search, variable neighborhood search, ant colony optimization, simulated annealing, genetic algorithms, scatter search, etc. Particularly appealing appear to be genetic algorithms, for which a number of different proposals were published, with special reference to how to optimize the crossover operator. As for the second line, the proposals are different, but they still have to settle and show how they compare on a broader range of problems, other than those for which they were originally presented. One example is the so-called *Forward and Backward (F&B)* approach [55], [56] which implements a memory-based look ahead strategy based on the past search history. The method iterates a partial exploration of the solution space by generating a sequence of enumerative trees of two types, called forward and backward trees, such that a partial solution of the forward tree has a bound on its completion cost derived from partial solutions of the backward tree, and vice-versa.

## III. Branch and price/column generation-based approaches

Also, Branch-and-price/column generation algorithms are usually adopted to solve set partitioning formulations. Branch-and-price algorithms have been proved to be successful for the exact solution of a wide variety of routing problems, including some of the most famous and classical ones, like the VRP and VRPTW, and are at the moment the exact leading methodology. While the branch-and-price scheme is a successful exact method and column generation is a building block of it, their use has been extended to obtain high performing and efficient heuristic algorithms. We call these heuristic approaches branch-and price/ column generation-based approaches. They have the common characteristic of using branch-and-price and/or column generation to build heuristic solutions. However, numerous schemes have been proposed in the literature which differ in terms of how columns are generated and/or of how they are used to obtain a feasible solution. [57],[58] classified this approach as 4 classes: *restricted master heuristics*, *heuristic branching approaches* and *relaxation-based approaches.*

- *restricted master heuristics*

One of the most used schemes in the class of branch-and-price/column generation based approaches is called restricted master heuristic. This scheme is typically embedded in a branch-and-price approach where the set partitioning formulation is solved on a subset of columns generated by the solution of the pricing problem, thus obtaining a feasible solution quickly. The restricted master heuristic is widely used in branch and-price approaches as it enables a quick improvement of bounds and thus a speed up of the exact solution procedure. Also, the can be used as heuristic algorithms to generate the columns. The column generation phase may be performed in two different ways: either a heuristic is used which does not take into account the dual information given by the solution of the restricted master problem, or the column generation algorithm is based on the dual information, but only a restricted set of columns is generated. Most of the approaches belong to the first class. This is due to the fact that these approaches are much easier to implement as they simply require a heuristic scheme to generate columns and the set partitioning model. We first analyze the approaches based on heuristic column generation and then describe those based on the use of the dual information provided by the master problem.

- *heuristic branching approaches*

Heuristic branching approaches are branch-and-price algorithms where, in order to speed up the convergence of the solution method, branching is performed heuristically with the aim of pruning a large number of nodes of the branch-and-bound tree and thus obtaining a good solution quickly.

In column generation approaches and branch-and-price algorithms, it is important to have fast algorithms available for repeatedly solving the pricing sub-problem, i.e. identifying a variable (column) with negative reduced costs. For many hard problems, however, this sub-problem is also hard. Fast heuristics are therefore sometimes used for approaching the pricing problem. Note that it is fine when pricing in a column with negative reduced costs even when it is not one with minimum reduced costs. However, at the end of column generation it is necessary to prove that no further column with negative reduced costs exists, i.e. the pricing problem must finally be solved exactly. Otherwise, no quality guarantees can be given for the final solution of the whole column generation or branch-and-price algorithm, and they must be considered to be heuristic methods only. Most heuristic approaches for solving pricing problems are relatively simple construction methods. More sophisticated metaheuristics have so far been used less frequently.

Also, almost any effective B&B approach depends on some heuristic for deriving a promising initial solution, whose objective value is used as original upper bound. Furthermore, and as already mentioned, heuristics are typically also applied to some or all sub-problems of the B&B tree in order to eventually obtain new incumbent solutions and corresponding improved upper bounds. In order to keep the B&B tree relatively small, good upper bounds are of crucial interest. Therefore, metaheuristics are often also applied for these purposes. However, when performing a relatively expensive metaheuristic at each node of a large B&B tree in a straight-forward, independent way, the additional computational effort often does not pay off. Different calls of the metaheuristic might perform more or less redundant searches in similar areas of the whole search space. A careful selection of the B&B tree nodes for which the metaheuristic is performed and how much effort is put into each call is therefore crucial. As an example, [59], [60] describes a chunking-based selection strategy to decide at each node of the B&B tree whether or not reactive tabu search is called. The chunking-based strategy measures a distance between the current node and nodes already explored by the metaheuristic in order to bias the selection toward distant points. Reported computational results indicate that adding the metaheuristic improves the B&B performance.

- *relaxation-based approaches*

An optimal solution for a relaxation of the original problem often indicates in which areas of the original problem's search space good or even optimal solutions might lie. Solutions to relaxations are therefore frequently exploited in (meta-) heuristics.

Sometimes an optimal solution to a relaxation can be repaired by a problem specific procedure in order to make it feasible for the original problem and to use it as promising starting point for a subsequent metaheuristic (or exact) search. Often, the linear programming (LP) relaxation is used for this purpose, and only a simple rounding scheme is needed. For example, [46] combine interior point methods and metaheuristics for solving the multidimensional knapsack problem (MKP). In a first step an interior point method is performed with early termination. By rounding and applying several different ascent heuristics, a population of different feasible candidate solutions is generated. This set of solutions is then used as initial population for a path-relinking/scatter search. Obtained results show that the presented combination is a promising research direction.

Beside initialization, optima of LP relaxations are often exploited for guiding local improvement or the repairing of infeasible candidate solutions. For example, in [61] the MKP is considered, and variables are sorted according to increasing LP-values. A greedy repair procedure considers the variables in this order and removes items from the knapsack until all constraints are fulfilled. In a greedy improvement procedure, items are considered in reverse order and included in the knapsack as long as no constraint is violated.

Another possibility of exploiting the optimal solution of an LP relaxation is more direct and restrictive: Some of the decision variables having integral values in the LP-optimum are fixed to these values, and the subsequent optimization only considers the remaining variables. Such approaches are sometimes also referred to as core methods, since the original problem is reduced and only its "hard core" is further processed. Obviously, the selection of the variables in the core is critical. Another example for exploiting the LP relaxation within metaheuristics is the hybrid tabu search algorithm from [48]. Here, the search space is reduced and partitioned via additional constraints fixing the total number of items to be packed. Bounds for these constraints are calculated by solving modified LP relaxations. For each remaining part of the search space, tabu search is independently applied, starting with a solution derived from the LP relaxation of the partial problem. The approach has further been improved in [49] by additional variable fixing.

Also other relaxations besides the LP relaxation are occasionally successfully exploited in conjunction with metaheuristics. The principal techniques for such combinations are similar.

The relaxation-based approaches are characterized by the fact that a feasible solution to the problem is generated from the information provided by the optimal solution of a relaxation of the master problem. Column generation is used to solve the relaxation. Once the relaxed solution is obtained, a heuristic procedure is used to generate a feasible solution to the problem.

Overall, Branch-and-price/column generation-based matheuristics are becoming more and more popular. This is due to the success of branch-and-price algorithms developed for the exact solution of routing problems. The scientific community has developed a deep knowledge of column generation approaches, and this knowledge is nowadays transferred also to the development of heuristic schemes. A further advantage of branch-and-price/column generation-based approaches is that they are flexible and easily adaptable to different problem characteristics. Most of the algorithms adopt the idea of using a set partitioning formulation and rely on heuristic schemes for the generation of columns.

## GENERAL DISCUSSION AND CONCLUSION

We have surveyed a multitude of examples where more powerful optimization systems were constructed by combining mathematical programming techniques and metaheuristics. Many very different ways exist for such hybridizations, and we have classified them into several major methodological categories and also we brought some examples of using them in vehicle routing problem. The probably most traditional approach is to use some metaheuristic for providing high-quality incumbents and bounds to a B&B-based exact method. On the other hand, quickly solved relaxations or the primal-dual relationship are often used for guiding or narrowing the search in metaheuristics. A relatively new and highly promising stream are those methods in which B&B is modified in some way in order to follow the spirit of local search based metaheuristics. A nowadays frequently and successfully applied approach is large neighborhood search by means of ILP techniques. When extending this concept towards searching the neighborhood defined by the common and disjoint properties of two or more parental solutions, we come to solution merging approaches. Furthermore highly promising hybrid approaches are those where metaheuristics are utilized within more complex branch-and-cut and branch and- price algorithms for cut separation and column generation, respectively. As noted, some approaches from the literature can be considered to fall into several of the methodological categories we have identified. Although a lot of experience already exists with such hybrid systems, it is usually still a tough question which algorithms and kinds of combinations are most promising for a new problem at hand.